\documentclass[11pt]{amsart}

\usepackage[T1]{fontenc}
\usepackage[ansinew]{inputenc}

\usepackage{euler}

\usepackage{amsmath}  
\usepackage{amsthm}   
\usepackage{latexsym}
\usepackage{amsfonts}
\usepackage{txfonts}
\usepackage{mathrsfs}

\usepackage{graphicx}
\usepackage{psfrag}

\numberwithin{figure}{section}

\numberwithin{equation}{section}

\begin{document}

\newtheorem{defn}{Definition}[section]
\newtheorem{teor}[defn]{Theorem}
\newtheorem{prop}[defn]{Proposition}
\newtheorem{exmp}[defn]{Example}
\newtheorem{lemm}[defn]{Lemma}
\newtheorem{seur}[defn]{Corollary}
\newtheorem{remark}[defn]{Remark}

\title[Existence of an absolute minimizer via Perron's method]
{Existence of an absolute minimizer via Perron's method}

\author{Vesa Julin}

\address{Department of Mathematics and Statistics,
P.O.Box 35, FIN-40014 University of Jyv\"askyl\"a, Finland}
\email{vesa.julin@jyu.fi}


\keywords{supremal functionals, absolute minimizer}
\subjclass[2000]{49J45, 49J99}
\date{\today}

\begin{abstract} 
 In this paper the existence of an absolute minimizer for a functional 
\[
F(u,\Omega) =  \underset{x \in \Omega}{ \text{ess sup}} \,  f (x, u(x), Du(x)) 
\] 
is proved by using Perron's method. The function is assumed to be quasiconvex and uniformly coercive. This completes the result by Champion, De Pascale and Prinari \cite{Thierry}.

\end{abstract}

\maketitle

\section{Introduction}
\vspace*{4mm} 

There has been increasing interest in Calculus of variations for \( L^{\infty}\) functionals in recent years. By Calculus of variations for \( L^{\infty}\) functionals we mean minimizing problems involving functionals of the form
\begin{equation}
\label{functional}
F(u, \Omega) = \underset{x \in \Omega}{ \text{ess sup}} \,  f (x, u(x), Du(x) )   
\end{equation}
where \( \Omega \subset \mathbb{R}^n\) is a bounded domain, \( f:\Omega \times \mathbb{R} \times \mathbb{R}^n \to \mathbb{R}  \) is measurable function and \( u \) is (locally) Lipschitz continuous in \( \Omega \). At this point we would like to mention the pioneering works by Aronsson in 1960's (\cite{Aronsson1} , \cite{Aronsson2} and \cite{Aronsson3}). 

One of the fundamental problems in the area is the existence of a so called absolute minimizer for the functional (\ref{functional}) with a given Dirichlet boundary data. That is to find a function $u \in W_g^{1,\infty}(\Omega)$ such that for every $V \subset \subset \Omega$ it holds
\begin{equation}
\label{minprob}
 F(v, V) \geq F(u, V ) \qquad \text{if} \,\, v \in W^{1,\infty}(V) \cap C(\overline{V}) \,\, \text{and} \,\,v = u \,\, \text{on} \, \partial V. 
\end{equation}
Here $g \in W^{1,\infty}(\Omega)$ is a given function and $W_g^{1,\infty}(\Omega)$ denotes the space of function $u-g \in W^{1,\infty}(\Omega) \cap C_0(\Omega)$. 

There are basically two ways to find an absolute minimizer for (\ref{functional}): the one is the  \(L^p \) approximation argument and the other Perron's method. Bhattacharya, DiBenedetto and Manfredi were the first ones to use \(L^p \) approximation in \cite{BhaDiBMan} where they proved the existence of an absolute minimizer for (\ref{functional}) in the special case \( f(x,s,p) = |p|\). The same method was later used by Barron, Jensen, Wang \cite{BarronJensenWang2} and more recently by Champion, De Pascale, Prinari \cite{Thierry} for much more general type of \( L^{\infty}\) functionals. Essentially what they proved is that whenever \(f \) is quasiconvex with respect to the last variable and uniformly coercive (see conditions \( (H1) \) and \( (H2) \) in the next section) \( L^{\infty}\) variational problem has an absolute minimizer. 

The use of Perron's method in \( L^{\infty}\) Calculus of variations dates back in 1960's, when Aronsson \cite{Aronsson3}  himself proved the existence of absolute minimizer in the special case \( f(x,s,p) = |p|\). Similar treatment was done by Juutinen \cite{Petri} and Milman \cite{Milman} in general metric spaces. Champion, De Pascale and Prinari \cite{Thierry} showed that this method also gives the existence of absolute minimizer for more general type of functionals. However, this result is not as general as in the case of \(L^p \) approximation done by the same authors, since in addition to the natural conditions \((H1 )\) and \( (H2) \) for \( f \) another assumption \( (H3) \) was needed. 

A natural question is wheather we will be able to get as strong existence result with Perron's method as we get by using \(L^p \) approximation? In this paper we prove that this is indeed the case. We use Perron's method to prove our main result, Theorem \ref{mainteor}, which states that an absolute minimizer exists if the integrand $f$ satisfies the natural conditions \((H1 )\) and \( (H2) \). The key is to define function classes which we call absolute superminimizers and absolute subminimizers. This gives an easy way to characterize the solution and the proof becomes rather straightforward.

\section{Preliminaries}
\vspace*{4mm} 

As we said in the introduction, the key of the proof is to use the following definition.

\begin{defn}
 \label{abs.min.}
A function \( u \in W_{loc}^{1,\infty}(\Omega)\) is an absolute superminimizer (subminimizer) of functional (\ref{functional}) if for all \( V \subset \subset \Omega \) and for \( v \in W^{1,\infty}(V) \cap C(\overline{V}) \) such that \( v > u \, ( v < u) \) in \( V \) and \( v = u \) on \( \partial V \) it holds
\[
  F(v, V) \geq F(u, V ).
\]
A function is an absolute minimizer of (\ref{functional}) if it is both absolute super- and subminimizer. 
\end{defn}

It is easy to see that \(u \in W_{loc}^{1,\infty}(\Omega) \) is an absolute minimizer of (\ref{functional}) if and only if for all \( V \subset \subset \Omega \) we have that
\[
 F(v, V) \geq F(u, V ) \qquad \text{if} \,\, v \in W^{1,\infty}(V) \cap C(\overline{V}) \,\, \text{and} \,\,v = u \,\, \text{on} \, \partial V. 
\]
Therefore our definition of absolute minimizer in Definition \ref{abs.min.} coincides with the one introduced earlier in (\ref{minprob}).

Function \(f\) in (\ref{functional}) is assumed to be measurable and to satisfy the following conditions:

\vspace*{4mm}

\begin{itemize}
 \item[(H1)] \emph{For a.e. \(x \in \Omega \) the map \( f(x, \cdot, \cdot )\) is lower semicontinuous on \( \mathbb{R} \times \mathbb{R}^n\) and for all \( (x,s) \in \Omega \times \mathbb{R} \) \( f(x, s, \cdot) \) is quasiconvex on \(\mathbb{R}^n \) i.e. for all \( p,q \in \mathbb{R}^n \) and \( 0 \leq t \leq 1 \) it holds
\[
 f(x,s, tp + (1-t)q) \leq \max \{ f(x,s,p) , f(x,s,q) \} 
\]
for all \( (x,s) \in \Omega \times \mathbb{R} \),}
\item[(H2)] \emph{For all \( c \in \mathbb{R} \) there is \( R \geq 0\) such that for every \( (x,s) \in \Omega \times \mathbb{R} \) it holds
\[
 \{ p \in \mathbb{R}^n  \mid f(x,s,p) \leq c \} \subset B(0, R).
\]}
 \end{itemize}

Condition \( (H2) \) is just uniform coerciveness. Condition \( (H1) \) guarantees that our functional has the right kind of  semicontinuity property as the following result states. The proof can be found in \cite{Thierry}.
\begin{teor}
 \label{teor.barron}
Let \( f : \Omega \times \mathbb{R} \times \mathbb{R}^n \rightarrow \mathbb{R} \) satisfy \((H1) \) and \( (H2) \). Then the functional (\ref{functional}) is sequentially lower semicontuous in \( W^{1,\infty}(\Omega) \) with respect to weak*-convergence i.e.
\[
 \underset{j \to \infty}{\text{lim inf}} \, F(u_j, \Omega) \geq F(u, \Omega)
\]
whenever \( u_j \to u\) weakly* in   \( W^{1,\infty}(\Omega) \).  
\end{teor}
The previous result is the key tool in finding a minimizer for functional (\ref{functional}) in $ W_g^{1,\infty}(\Omega)$ by the direct method of Calculus of variations. We sketch the proof for readers convenience.
\begin{teor}
 \label{existence1}
Suppose \( f : \Omega \times \mathbb{R} \times \mathbb{R}^n \rightarrow \mathbb{R} \) satisfies \((H1) \) and \( (H2) \) and \(g \in W^{1,\infty}(\Omega) \). Then the functional (\ref{functional}) has at least one minimizer in \( W_g^{1,\infty}(\Omega) \) i.e. there is \( u \in W_g^{1,\infty}(\Omega)\) such that 
\[
 F(u, \Omega) = \inf \, \{ F(v, \Omega) \mid v \in W_g^{1,\infty}(\Omega) \} .
\]

\end{teor}

\textit{Proof.} Denote \( c = F(g, \Omega) \) and let \( (u_j) \) be a sequence such that 
\[
 F(u_j, \Omega) \to \inf_{v \in W_g^{1,\infty}(\Omega)} F(v, \Omega) \in [- \infty, c].
 \]
Obviously the sequence can be chosen so that \( F(u_j, \Omega) \leq c \) for all \( j \) and therefore the condition \((H2) \) implies that \( (u_j) \) is bounded in \(W^{1,\infty}(\Omega) \). Hence \( \inf_{v \in W_g^{1,\infty}(\Omega)} F(v, \Omega) > - \infty \). Moreover we may assume that the sequence \(u_j\) weakly*-converges towards some \( u \in W_g^{1,\infty}(\Omega)\). Theorem \ref{teor.barron} guarantees that the function \( u \) is a minimizer of \( F( \cdot, \Omega) . \qquad \Box \)

\smallskip

Champion, De Pascale and Prinari (\cite{Thierry},  \emph{ Theorem 4.1}) proved that the conditions \((H1) \) and \((H2) \)  pretty much guarantees the existence of an absolute minimizer. To be quite precise, they need yet to assume that the integrand $f$ is continuous with respect to its second variable. This is the proof which uses \( L^p\) approximation. But when they use Perron's method  (\cite{Thierry}, \emph{Theorem 4.7}) the following additional assumption is needed: 
\begin{itemize}
 \item[(H3)] \emph{For any open subset \( V \subset \Omega \,, g \in W^{1,\infty}(V) \cap C(\overline{V}) \) and \( y \in V\), the image set  
\[ 
\{ u(y) \mid  u \, \, \text{is minimizer of} \, F( \cdot, V) \,  \text{in} \, W_g^{1,\infty}(V) \} \subset \mathbb{R}
\] 
is connected.}
\end{itemize}
The point of this paper is to prove the existence of an absolute minimizer by using Perron's method, whitout using the assumption \( (H3) \). 
\begin{teor}
 \label{mainteor}
Suppose \( f : \Omega \times \mathbb{R} \times \mathbb{R}^n \rightarrow \mathbb{R} \) satisfies \((H1) \) and \((H2) \). Then for any \( g \in W^{1,\infty}(\Omega)\) the functional (\ref{functional}) has at least one absolute minimizer in \(W_g^{1,\infty}(\Omega) \).
\end{teor}

\section{Existence of absolute minimizer}
\vspace*{4mm} 

The outline of the proof of Theorem \ref{mainteor} is quite standard. We will construct our absolute minimizer piece by piece by using Theorems \ref{teor.barron} and \ref{existence1}.

We will frequently use the following notations,
\[
 \begin{split}
  \mathscr{A}(g,\Omega) &= \{ u \in W_g^{1,\infty}(\Omega) \mid u \, \text{is a minimizer of} \, F(\cdot, \Omega) \, \text{in} \,  W^{1,\infty}(\Omega)  \} \\
  \mathscr{A}_{sup}(g, \Omega) &= \mathscr{A}(g,\Omega) \cap \{ u \in W_g^{1,\infty}(\Omega) \mid u \, \text{is an absolute superminimizer of} \, F(\cdot, \Omega) \} \\
\mathscr{A}_{sub}(g, \Omega) &= \mathscr{A}(g,\Omega) \cap \{ u \in W_g^{1,\infty}(\Omega) \mid u \, \text{is an absolute subminimizer of} \, F(\cdot, \Omega) \}.
 \end{split}
\]

Since we are using Perron's method, it is rather obvious that the following two lemmas are needed.

\begin{lemm}
 \label{lemma2}
Assume that  \(g \in W^{1,\infty}(\Omega) \).
\begin{itemize}
 \item [(i)] Suppose that \( u_1, u_2 \in \mathscr{A}_{sub}(g, \Omega) \) and set \( u = \max(u_1, u_2) \). Then \( u \in \mathscr{A}_{sub}(g, \Omega)\).
 \item [(ii)] Suppose that \( u_1, u_2 \in \mathscr{A}_{sup}(g, \Omega) \) and set \( u = \min(u_1, u_2) \). Then \( u \in \mathscr{A}_{sup}(g, \Omega)\).
\end{itemize}
\end{lemm}

\textit{Proof.} \( (i) \) First of all, since \( u_1, u_2 \in \mathscr{A} (g, \Omega) \) we have that \( u \in \mathscr{A} (g, \Omega) \). Indeed, denote \\ 
\( W = \{ x \in \Omega \mid  u_1(x) > u_2(x) \}  \) and deduce that
\[
 \begin{split}
  F( u , \Omega ) &=  \max \{ F( u_1 , W ), \, F( u_2, \Omega \backslash W ) \} \leq  \max \{ F( u_1 , \Omega ), \, F( u_2, \Omega ) \} \\
		&= \inf \, \{ F(v, \Omega) \mid v \in W_g^{1,\infty}(\Omega) \} .
 \end{split}
\]
Hence \( u \in \mathscr{A} (g, \Omega) \).

Suppose next that \( V \subset \subset \Omega \) and \( h \in W^{1,\infty}(V) \cap C(\overline{V})\) is such that \( h < u \) in \( V \) and \( h = u \) on \( \partial V \). Divide \( V \) into two parts in two ways. Set first
\[
 \begin{split}
  U_1 &= \{ x \in V \mid h(x) < u_1(x) \} \\
  U_2 &= \{ x \in V \mid h(x) < u_2(x) \}.
 \end{split}
\]
Notice that \( h = u_1 \) on \( \partial U_1 \) and \( h = u_2 \) on \( \partial U_2 \). Define next
\[
 \begin{split}
  W_1 &= \{ x \in V \mid  u_1(x) > u_2(x) \} \\
  W_2 &= \{ x \in V \mid  u_1(x) \leq u_2(x) \}.
 \end{split}
\]
It is immediate that  \( u = u_1 \) in \( W_1 \), \( u = u_2\) in \( W_2\) and \( W_1 \subset U_1, W_2 \subset U_2\). Together with the fact that $u_1$ and $u_2$ are absolutely subminimizers these imply
\[
\begin{split}
 F( h, V ) 	&= \max (F( h, U_1 ), F( h, U_2 )) \geq \max (F( u_1, U_1 ), F( u_2, U_2 )) \\
		&\geq \max (F( u_1, W_1 ), F( u_2, W_2 )) = \max (F( u, W_1 ), F( u, W_2 )) = F( u, V ).	
 \end{split}
\]
  
Part \( (ii) \) goes similarly. \( \qquad \Box\)

\begin{lemm}
\label{lemma1}
Fix  \(g \in W^{1,\infty}(\Omega) \). Then the set \(\mathscr{A}(g,\Omega) \) is non-empty. Moreover, consider the functions 
\[
 w(x) \coloneqq \sup_{v \in \mathscr{A}(g,\Omega)} v(x)
\]
and 
\[
u(x) \coloneqq \inf_{v \in \mathscr{A}(g,\Omega)} v(x).
\]
Then \( w \in \mathscr{A}_{sup}(g,\Omega)\) and \(u \in \mathscr{A}_{sub}(g,\Omega)\).
\end{lemm}

\textit{Proof.} The fact that \(\mathscr{A}(g,\Omega) \) non-empty is just theorem \ref{existence1}. We will only prove that \( w \in \mathscr{A}_{sup}(g,\Omega) \), since the proof for \(u \) is completely analogous. 

Suppose \( \{x_1, x_2, \dots \} \) is dense in \( \Omega\). Set \( k \in \mathbb{N} \) and find functions \( u_1^k, \dots, u_k^k \in \mathscr{A}(g,\Omega) \) such that \( u_i^k(x_i) + \frac{1}{k} \geq w(x_i) \geq u_i^k(x_i) \) for every \( i = 1, \dots,k\). Denote \( v_k(x) = \max ( u_1^k(x), \dots, u_k^k(x) ) \). Looking at the first part of the proof of Lemma \ref{lemma2} we conclude that \( v_k \in \mathscr{A}(g,\Omega) \). By doing this for all \( k \) we obtain a sequence \( (v_k)\) of minimizers of (\ref{functional}) such that \(v_k(x) \to w(x) \) pointwise in a dense subset of \( \Omega \). Since functions \( v_k \) are minimizers we have \( F( v_k , \Omega) \leq F( g , \Omega) \) for all \( k \). By \( (H2) \) the sequence is bounded in \( W^{1,\infty}(\Omega) \) and by passing to a subsequence we may assume that \( v_k \overset{w*}{\to} w \). Since \( F( \cdot, \Omega) \) is weakly* lower semicontinuous we conclude that \(w \in \mathscr{A}(g,\Omega) \) . 

Suppose that \( w\) is not an absolute superminimizer. Then there would be \( V \subset \subset \Omega \) and \( h \in W^{1,\infty}(V) \cap C(\overline{V}) \) such that \( h > w \) in \( V  \, , h = w \) on  \( \partial V \) and
\[
 F(h, V) < F(w, V) .
\]
Define
\[
 \tilde{w}(x) = \begin{cases}
                 h(x) & x \in V \\
		 w(x) & x \in \Omega \backslash V .
                \end{cases}
\]
Then \( \tilde{w} \in W_g^{1,\infty}(\Omega)\) and 
\[
 F(\tilde{w}, \Omega) = \max (F(\tilde{w}, V), F(\tilde{w}, \Omega \backslash V)) =  \max (F( h, V), F( w , \Omega \backslash V)) \leq F(w, \Omega) .
\]
Hence \(\tilde{w} \in \mathscr{A}(g,\Omega)\) and \( \tilde{w} > w \) in \( V \). But this contradicts the definition of \( w \), and \( w \) is therefore an absolute superminimizer . \( \qquad \Box \)

Now we are ready to prove our main result.

\vspace*{4mm}

\textbf{Proof of Theorem \ref{mainteor}:} \, Fix \(g \in W^{1,\infty}(\Omega)\). We will show that absolute minimizer can be found by the formula
\[
 \bar{u}(x) = \inf_{v \in \mathscr{A}_{sup}(g,\Omega)} v(x)
\]
or
\[
 \bar{w}(x) = \sup_{v \in \mathscr{A}_{sub}(g,\Omega)} v(x).
\]
Moreover \( \bar{u}\) is the smallest and \( \bar{w} \) is the biggest absolute minimizer of functional (\ref{functional}) in \(  W_g^{1,\infty}(\Omega)\). We will only show that \(\bar{u} \) is an absolute minimizer, since the proof for \(\bar{w} \) is completely analogous.

\vspace*{4mm} 

\textbf{Claim 1:} \, \( \bar{u} \in \mathscr{A}_{sup}(g,\Omega) \).

Just like in the proof of Lemma \ref{lemma1} we choose a dense subset \( \{x_1, x_2, \dots \} \) of \( \Omega \) and for all \( k \in \mathbb{N} \) functions \( u_1^k, \dots , u_k^k \in \mathscr{A}_{sup} \) such that \( u_j^k(x_j) \geq \bar{u}(x_j) \geq u_j^k(x_j) - \frac{1}{k} \) for \( j = 1, \dots , k \). Denote \( v_k(x) = \min (u_1^k(x), \dots , u_k^k(x) )\). Then by Lemma \ref{lemma2} \( v_k \in \mathscr{A}_{sup}(g,\Omega) \) and by construction \( v_k \to \bar{u} \) pointwise in a dense subset of \( \Omega \). In particular, \( v_k \in \mathscr{A}(g,\Omega) \) and therefore \( F( v_k , \Omega) \leq F( g , \Omega) \) for all \( k \). Again by \( (H2) \) the sequence is bounded in \( W^{1,\infty}(\Omega) \) and we may assume that \[ 
v_k \overset{w*}{\to} \bar{u} .
\] 
In particular, \( v_k \to \bar{u} \) uniformly. Moreover we may assume that the sequence is nonincreasing by considering \( \tilde{v}_k(x) = \min (v_1(x), \dots , v_k(x) )\).  

Suppose \( V \subset \subset \Omega \) and \( h \in W^{1,\infty}(V) \cap C(\overline{V})\) is such that \( h > \bar{u} \) in \( V \) and \( h = \bar{u} \) on \(\partial V \). Denote 
\[
 V_k = \{ x \mid h(x) > v_k(x) \} .
\]
Since \( (v_k) \) is nonincreasing and converges uniformly to \( \bar{u} \) we have \(V_k \subset V_{k + 1} \subset \dots \subset V \) for all \( k \in \mathbb{N} \) and \( V_k \) is non-empty when \( k \) is large. Therefore 
\[
 F(h,V_k) \geq F(v_k, V_k)
\]
since \(v_k \in \mathscr{A}_{sup}(g,\Omega) \). Fix a large \(  k_0 \) for a moment. For all \( k \geq  k_0 \) we have 
\[
 F(h,V) \geq F(h, V_k) \geq F(v_k, V_k) \geq F(v_k, V_{k_0}).
\]
Therefore letting \( k \to \infty \) we have by the weak* semicontinuity of \( F \) that
\[
 F(h,V) \geq {\lim \inf}_{k \to \infty} F(v_k, V_{k_0}) \geq F(\bar{u}, V_{k_0}) .
\]
Finally by letting \( k_0 \to \infty\) we conclude
\[
 F(h,V) \geq F(\bar{u}, V)
\]
which implies \( \bar{u} \in \mathscr{A}_{sup}(g,\Omega) \).

\vspace*{4mm} 

\textbf{Claim 2:} \, \( \bar{u} \in \mathscr{A}_{sub}(g,\Omega) \).

We prove this by contradiction. Suppose that there is \( V' \subset \subset \Omega \) and  \( u \in  W^{1,\infty}(V') \cap C(\overline{V'})\) such that \( u < \bar{u} \) in  \( V'\), \( u = \bar{u} \) on \( \partial V' \) and 
\begin{equation}
\label{theor1}
F(u,V') < F( \bar{u} ,V') .
\end{equation}
Define
\[
 w(x) = \sup_{v \in \mathscr{A}( \bar{u}, V')} v(x) .
\]
By Lemma \ref{lemma1} \(w \in   \mathscr{A}_{sup}(\bar{u}, V') \). In particular \( w \) is a minimizer of functional \( F( \cdot, V') \) and therefore by Claim 1 and (\ref{theor1}) the set
\begin{equation}
\label{definitionV}
 V = \{ x \in V' \mid w(x) < \bar{u}(x) \}
\end{equation}
is non-empty. Define a function
\begin{equation}
\label{definition.uhat}
 \hat{u}(x) = \begin{cases}
                 w(x) & x \in V \\
		 \bar{u}(x) & x \in \Omega \backslash V .
                \end{cases}
\end{equation}
Our goal is to show that \( \hat{u} \in \mathscr{A}_{sup}(g, \Omega) \), which contradicts the definition of \( \bar{u} \) and the claim will then follow.
 
Therefore assume that \( U \subset \subset \Omega \), \( h > \hat{u} \) in \( U \) and \( h = \hat{u} \) on \( \partial U \). Denote
\begin{equation}
\label{definitionW}
 W = \{ x \in U \mid h(x) > \bar{u}(x) \}.
\end{equation}
Since \( \bar{u} \in \mathscr{A}_{sup}(g,\Omega) \) we have 
\begin{equation}
 \label{theor2}
\begin{split}
 F(h,U) &= \max \{ F(h, W), F(h, U \backslash W)\} \geq \max \{ F( \bar{u}, W), \, F(h, U \backslash W)\} \\
	&= \max \{ \underbrace{F( \bar{u}, W \backslash V)}_{*}, \, \underbrace{F( \bar{u}, W \cap V), F(h, U \backslash W)}_{**}\}.
\end{split}
\end{equation}

Consider first the term \( (* )\). Definition (\ref{definitionW}) yields \( W \subset U \) and thereby  \( W \backslash V \subset U \backslash V\). Suppose that \( x \in U \backslash V \). By (\ref{definition.uhat}) we have \( \hat{u}(x) = \bar{u}(x)\) which implies \( h(x) > \hat{u}(x) = \bar{u}(x)\). This implies \( x \in W \) and therefore \(  U \backslash V \subset W \backslash V \). Hence \( W \backslash V = U \backslash V \) and in particular
\begin{equation}
\label{theor3}
  F( \bar{u}, W \backslash V) = F( \bar{u}, U \backslash V) .
\end{equation}

Next consider the term \( (**)\). Using the definitions of (\ref{definitionV}) and (\ref{definitionW}) it is easy to see that
\begin{equation}
 \label{theor4}
\begin{split}
U \cap V \cap W &= V \cap W \,\, \text{and} \\
(U \cap V ) \backslash W &= U \backslash W .
\end{split}
\end{equation}

Next we notice that \( \min(\bar{u}(x), h(x)) > w(x) \) for \( x \in U \cap V \) and \(\min(\bar{u}(x), h(x)) = w(x) \) for \(  x \in \partial (U \cap V) \). Since \( w \) is an absolute superminimizer in \( V \) we have
\begin{equation}
 \label{theor5}
\begin{split}
 F(w, U \cap V) &\leq F( \min(\bar{u}, h) , U \cap V) = \max \{ F( \bar{u} , (U \cap V \cap W)) , \, F( h , (U \cap V)  \backslash W ) \} \\
		&=  \max \{ F( \bar{u} ,  V \cap W ), \, F( h , U  \backslash W ) \} 
\end{split}
\end{equation}
where the last equality follows from (\ref{theor4}).

Combining (\ref{theor2}), (\ref{theor3}) and (\ref{theor5}) yields
\[
 F(h, U) \geq \max \{ F( \bar{u}, U \backslash V), \, F(w, U \cap V)\} = F( \hat{u}, U )
\]
 since \( \hat{u}(x) = \bar{u}(x)  \) for \( x \in  U \backslash V \) and \(\hat{u}(x) = w(x) \) for \( x \in V . \qquad \Box\)

\smallskip

\textbf{Acknowledgments}. The author wishes to thank the referee for his comments on the article. The author was partially supported by the Academy of Finland, project \#129784.


\begin{thebibliography}{BarronJensenWang2}

\bibitem{Aronsson1} Aronsson G., \emph{Minimization problems for the functional \( \sup_x F(x,f(x), f'(x)) \)}, Arkiv för Mat. \textbf{6} (1965), 33--53.

\bibitem{Aronsson2} Aronsson G., \emph{Minimization problems for the functional \( \sup_x F(x,f(x), f'(x)) \)} part 2, Arkiv för Mat. \textbf{6} (1966), 409--431.

\bibitem{Aronsson3} G. Aronsson, \emph{Extension of functions satisfying Lipschitz conditions}, Arkiv för Mat. \textbf{6} (1967), no. 28, 551--561.

\bibitem{BarronJensenWang2} Barron E.N. , Jensen R. , Wang C.Y. , \emph{The Euler equation and absolute minimizers of \( L^{\infty}\) functionals}, Arch. Ration. Mech. Anal. \textbf{157} (2001), no. 4, 255--283.

\bibitem{BhaDiBMan} Bhattacharya T. , DiBenedetto E. , Manfredi J. , \emph{Limits as \(p \to \infty \) of \(\triangle_p u_p = f  \) and relatedextremal problems} Rendiconti del Sem. Mat. Fascicolo Speciale Non Linear PDE's, Univ. di Torino, (1989), 15-68.

\bibitem{Thierry} Champion T. , De Pascale L. , Prinari F. , \emph{\( \Gamma \)-convergence and absolute minimizers for supremal functionals}  ESAIM Control Optim. Calc. Var. \textbf{10}   (2004),  no. 1, 14--27 (electronic).

\bibitem{Petri} Juutinen P. , \emph{Absolutely minimizing Lipschitz extensions on a metric space},  Ann. Acad. Sci. Fenn. Math.  \textbf{27}  (2002),  no. 1, 57--67.

\bibitem{Milman}  Milman V. A., \emph{Absolutely minimal extensions of functions on metric spaces}, (Russian) Mat. Sb., \textbf{190} (1999), no 6, 83--110;  translation in Sb. Math.  \textbf{190}  (1999),  no. 5-6, 859--885.  

\end{thebibliography}
\end{document}